\documentclass[12pt]{amsart}
 
\usepackage[dvips]{graphicx}
	\usepackage{amsmath,amssymb,amsthm,wrapfig,amsfonts,enumerate,citesort}

	\newtheorem{dfn}{Definition}[section]
	\newtheorem{thm}[dfn]{Theorem}

	\newtheorem{lem}[dfn]{Lemma}
	\newtheorem{rem}[dfn]{Remark}

        \newtheorem{bai}{Case}
	\newtheorem{ack}{Acknowledgements\!\!}

	\numberwithin{equation}{section}

	\def\notin{\not\in}

	\newcommand{\dist}{\mathop{\mathit{d}} \nolimits}
	
	\newcommand{\diam}{\mathop{\mathrm{diam}} \nolimits}

        \newcommand{\e}{\mathop{\varepsilon}              \nolimits}

\begin{document}

	\title[Two infinite versions of nonlinear Dvoretzky's theorem]
    {Two infinite versions of nonlinear Dvoretzky's theorem}
	\author[Kei Funano]{Kei Funano}
	\address{Research Institute for Mathematical Sciences, Kyoto University,
Kyoto 606-8502 JAPAN}
	\email{kfunano@kurims.kyoto-u.ac.jp}
	\thanks{This work was partially supported by Grant-in-Aid for
    Research Activity (startup) No.23840020.}
	\dedicatory{}
	\date{\today}

	\maketitle


\begin{abstract}We make two additions to recent results of Mendel and
 Naor on nonlinear versions of Dvoretzky's theorem. We consider the cases of metric spaces with
 infinite Hausdorff dimension and countably infinite metric spaces. 
 \end{abstract}
	\setlength{\baselineskip}{5mm}

\section{Introduction and the statement of the results}
We say that a metric space $X$ is \emph{embedded with distortion} $D\geq 1$
in a metric space $Y$ if there exist a map $f:X\to Y$ and a constant $r>0$ such that
\begin{align*}
 r \dist_X(x,y) \leq \dist_Y(f(x),f(y))\leq D r \dist_X(x,y) \text{ for all
 }x,y\in X.
 \end{align*}Such a map $f$ is called a \emph{$D$-embbeding}.

 Dvoretzky's theorem states that for every $\e>0$, every
 $n$-dimensional normed space contains a $k(n,\e)$-dimensional subspace
 that embeds into a Hilbert space with distortion $1+\e$
 (\cite{Dvo60}). This theorem was conjectured by Grothendieck
 (\cite{groth}). See \cite{Mil71,Mil92,MilSch96,Sch2,Sch} for the estimates of $k(n,\e)$
 and the further developments related to this theorem. 

  Bourgain, Figiel, and Milman proved the following theorem as a natural
  nonlinear variant of Dvoretzky's theorem.
 
 \begin{thm}[{cf.~\cite{BFM86}}]\label{bfmt}There exists two universal constants
  $c_1,c_2>0$ satisfying the following. For every $\e>0$ every finite
  metric space $X$ contains a subset $S$ which embeds into a Hilbert
  space with distortion $1+\e$ and
  \begin{align*}
   |S|\geq \frac{c_1 \e}{\log (c_2/\e)} \log |X|.
   \end{align*}
  \end{thm}

  See \cite{blmn}, \cite{MN07}, \cite{NT} for the further investigation.
  It is natural to try to get some versions of the above theorem in the case where $|X|=\infty$. In this paper we prove the following. 
\begin{thm}\label{finiteTHM}
 For every $\e>0$, every countable infinite metric space $X$ has an
 infinite subset which embeds into an ultrametric space with distortion
 $1+\varepsilon$.
 \end{thm}Recall that a metric space $(U,\rho)$ is called an
 \emph{ultrametric space} if for every $x,y,z\in X$ we have $\rho(x,y)
 \leq \max\{ \rho (x,z), \rho(z,y) \}$. Since every separable ultrametric space isometrically
 embeds into a Hilbert space (\cite{VeTi}), we verify that Theorem
 \ref{bfmt} holds in the case where $|X|=\infty$.

 Recently Mendel and Naor proved another variant of Dvoretzky's
 theorem, answering a question by T.~Tao (\cite{MN11}). For a metric
 space $X$ we denote its Hausdorff dimension by $\dim_H(X)$. A subset of
 a complete separable metric space is called an \emph{analytic
 set} if it is an image of a complete separable metric space under a
 continuous map. Note that analytic sets are not necessarily
 complete. For example any Borel subsets of a complete separable metric space
 are analytic sets (refer to \cite{kechris} for analytic sets). 
 \begin{thm}[{cf.~\cite[Theorem 1.7]{MN11}}]\label{MNDV}There exists a universal constant $c\in
  (0,\infty)$ such that for every $\e \in (0,\infty)$, every analytic
  set $X$ whose Hausdorff
  dimension is finite has a closed subset $S\subseteq X$ that embeds with distortion
  $2+\e$ in an ultrametric space, and
  \begin{align*}
   \dim_H(S)\geq \frac{c\e}{\log (1/\e)}\dim_H(X).
   \end{align*}
  \end{thm}

  In \cite{MN11} Mendel and Naor stated the above theorem only for
  compact metric spaces. As they remarked in \cite[Introduction]{MN11},
  their theorem is valid for more general metric spaces. For example
  the above theorem holds for every analytic set $X$ since the problem
  can be reduced to the case of a compact
  subset of $X$ with the same Hausdorff dimension (see \cite{carleson},
  \cite[Corollary 7]{how}). 

  In the following theorem we consider the case where
  $\dim_H(X)=\infty$.

 \begin{thm}\label{THM}For every $\e \in (0,\infty)$, every analytic
  set $X$ whose Hausdorff
  dimension is infinite has a closed subset $S$ such that $S$ can be embedded
  into an ultrametric space with distortion $2+\varepsilon$ and has infinite
  Hausdorff dimension. 
  \end{thm}

  It follows from the proof of Theorem \ref{MNDV} in \cite{MN11} that if $\dim_H(X)=\infty$, then $X$ contains an arbitrary large dimensional
  closed subset which embeds into an ultrametric space. Combining
  Theorem \ref{MNDV} with Theorem \ref{THM} we find that nonlinear
  Dvoretzky's theorem holds for all analytic sets. 

  The following theorem due to Mendel and Naor asserts that we cannot replace the distortion
  strictly less than $2$ in Theorems \ref{MNDV} and \ref{THM}.
\begin{thm}[{cf.~\cite[Theorem 1.8]{MN11}}]\label{hanrei}For every $\alpha>0$ there exists a compact metric space $(X,\dist)$ of Hausdorff
dimension $\alpha$, such that if $S\subseteq X$ embeds into a Hilbert space with distortion strictly smaller
than $2$ then $\dim_H(S) = 0$.
\end{thm}

  Theorem \ref{hanrei} immediately implies that the same result holds in
  the case where $\alpha=\infty$.

  It is known that $\ell_2$ does not embed into $\ell_p$ with finite
  distortion for any $p\in [1,\infty)\setminus \{2\}$ (\cite[Corollary 2.1.6]{alkal}). In particular, an infinite dimensional
  analogue of Dvoretzky's theorem is no longer true in the linear
  setting. In contrast to this fact, Theorem \ref{THM} asserts that an
  infinite dimensional Dvoretzky's theorem holds in the nonlinear
  setting.

  \section{Proof}

  We need the following lemma.
  \begin{lem}\label{lem1}Let $X$ be a separable metric space such that
   $\dim_H(X)=\infty$. Then there exists a sequence $\{
   K_i\}_{i=1}^{\infty}$ of mutually disjoint closed subsets
   of $X$ such that $\lim_{i\to
   \infty}\diam K_i=0$ and $\lim_{i\to \infty}\dim_H K_i =\infty$.
   \begin{proof}
   For every $x\in X$ we take a closed neighborhood $K_x$ of $x$
   such that $\diam K_x \leq 1$. Since $X$ is separable, applying the Lindel\"of covering
   theorem we get a countable subset
   $F\subseteq X$ such that $X= \bigcup_{x\in F}K_x$. Since
   $\dim_H(\bigcup_{x\in F} K_x)=\sup_{x\in F} \dim_H K_x $, there
    exists $x_1\in F$ such that $\dim_H K_{x_1}=\infty$ or there exists a sequence $\{ y_i \}_{i=1}^{\infty}\subseteq F$ such that
   $\{\dim_H K_{y_i}\}_{i=1}^{\infty}$ is strictly increasing and
    $\lim_{i\to \infty}\dim_H K_{y_i}=\infty$.

    We first consider the latter case. We put $K_1:=K_{y_1}$. By the
    monotonicity of $\dim_H K_{y_i}$ we have
    $\dim_H(K_{y_i}\setminus \bigcup_{j=1}^{i-1} K_{y_j})=\dim_H
    K_{y_i}$ for $i\geq 2$. Covering $K_{y_i}\setminus
    \bigcup_{j=1}^{i-1} K_{ y_j}$ by countably many closed subsets of
    diameter $\leq 1/i$, we thus find a closed subset $K_i\subseteq K_{y_i}\setminus
    \bigcup_{j=1}^{i-1} K_{y_j}$ such that $\dim_H K_i=\dim_H K_{y_i}$ and $\diam K_i \leq
    1/i$. This $\{K_i\}_{i=1}^{\infty}$ is a desired sequence. 

    We consider the former case. Covering  $K_{x_1}$ by countably many closed subsets $\{ K_{y}^1
    \}_{y\in F_1}$ so that $\diam K_y^1\leq 2^{-1}\diam K_{x_1}$, we
    have the following two cases: There exists $x_2\in F_1$ such that
    $\dim_H(K_{x_2}^1)=\infty $ or there exists a sequence $\{ y_i \}_{i=1}^{\infty}\subseteq F_1$ such that
   $\{\dim_H K_{y_i}^1\}_{i=1}^{\infty}$ is strictly increasing and
    $\lim_{i\to \infty}\dim_H K_{y_i}^1=\infty$. Since we have already proved the
    lemma in the latter case, we consider the former case. Continuing
    this process we may assume that there exists a chain
    $K_{x_2}^1\supseteq K_{x_3}^2 \supseteq K_{x_4}^3 \supseteq \cdots$
    of closed subsets of $X$ such
    that $\dim_H(K_{x_i}^{i-1})=\infty$ and $\diam K_{x_{i+1}}^{i}\leq
    2^{-1}\diam K_{x_i}^{i-1}$. Since $K_{x_i}^{i-1}\setminus \bigcup_{j=i}^{\infty}
     (K_{x_j}^{j-1}\setminus K_{x_{j+1}}^j)$ consists of at most one point, we get $\limsup_{i\to \infty} \dim_H(K_{x_i}^{i-1} \setminus
    K_{x_{i+1}}^{i})=\infty$. By taking a subsequence we may assume that
    $\lim_{i\to \infty} \dim_H(K_{x_i}^{i-1} \setminus
    K_{x_{i+1}}^i)=\infty$. Taking a closed subset $K_i\subseteq
    K_{x_i}^{i-1}\setminus K_{x_{i+1}}^i $ such that $\dim_H K_i \geq 2^{-1}
    \dim_H (K_{x_i}^{i-1}\setminus K_{x_{i+1}}^i)$ we easily see that
    this $\{ K_i\}_{i=1}^{\infty}$ is a desired sequence. This completes
    the proof.
    \end{proof}
   \end{lem}
 
 We first prove Theorem \ref{THM}. It turns out that Theorem \ref{finiteTHM} follows from
 the proof of Theorem \ref{THM}. 
  \begin{proof}[Proof of Theorem \ref{THM}]We take a
   sequence $\{ K_i\}_{i=1}^{\infty}$ of closed subsets
   of $X$ in Lemma \ref{lem1}. For each $i$ we fix an
   element $x_i \in K_i$. Note that closed subsets of analytic sets are
   also analytic sets. According to Theorem \ref{MNDV} there exist $A_i\subseteq K_{i}$ such that $\lim_{i\to \infty}\dim_H A_i =
   \infty$ and $A_i$ embeds into some ultrametric space $(U_i,\rho_i)$
   with distortion $2+\e$, i.e., there exist $f_i:A_i \to U_i$ satisfying
   \begin{align}\label{ms1}
    \dist(x,y)\leq \rho_i(f_i(x),f_i(y)) \leq (2+\e) \dist (x,y) \text{ for
    any }x,y \in A_i.
    \end{align}

   We divide the proof into three cases.

   \begin{bai}\label{case1}
    $ \{x_i \}_{i=1}^{\infty}$ is not bounded. 
    \end{bai}
   By taking a subsequence we may assume that $\lim_{n\to \infty}\dist
   (x_1,x_i)=\infty$ and $\diam K_{i} \leq 1/(2+\e)$. By taking a
   subsequence we may also assume that
   \begin{align}\label{ms2}
    1\leq \min\Big\{  \frac{\sqrt{1+\e}-1}{\sqrt{1+\e}
    \sqrt{1+2^{-1}\e}} ,  \frac{\sqrt{1+\e}-\sqrt{1+2^{-1}\e }}{\sqrt{1+2^{-1}\e}} , \frac{\sqrt{1+2^{-1}\e}-1}{2}         \Big\} \dist(A_1,A_2)
    \end{align}and
   \begin{align}\label{ms3}
    \dist(A_1,A_{i-1})\leq \min\Big\{  \frac{\sqrt{1+\e}-1}{\sqrt{1+\e}
    \sqrt{1+2^{-1}\e}} ,  \frac{\sqrt{1+\e}-\sqrt{1+2^{-1}\e }}{\sqrt{1+2^{-1}\e}}           \Big\} \dist(A_1,A_i).
    \end{align}for any $i\geq 2$. Put $R_i:=\dist(A_i,A_1)$ for $i\geq 2$. Note that $\diam f_i(A_i)\leq 1$ since $f_i$ satisfies (\ref{ms1}) and
   $\diam A_i \leq \diam K_{i} \leq 1/(2+\e)$.

   For each $i\ge 2$ we add an additional point $u_{i,0}$ to $f_i(A_i)$ and
   put $Y_i:=f_i(A_i)\cup \{ u_{i,0} \}$. Define the distance function
   $\tilde{\rho}_i $ on $Y_i$ as follows: $\tilde{\rho}_i(u,u_{i,0}):=R_i$ for
   $u\in f_i(A_i)$ and $\tilde{\rho}_i(u,v):= \rho_i(u,v)$ for $u,v\in
   f_i(A_i)$. Since $\diam f_i(A_i)\leq 1\leq R_i$, each
   $(Y_i,\tilde{\rho}_i)$ is an ultrametric space. Let us consider the space
   \begin{align}\label{ultrasp}
    U:= \{ (u_i) \in \prod_{i= 2}^{\infty} Y_i \mid  u_i\neq
    u_{i,0}\text{ only for finitely many }i         \}
    \end{align}and define the distance function $\rho$ on $U$ by
   \begin{align}\label{ultramet}
    \rho((u_i),(v_i)):= \sup_i \tilde{\rho}_i(u_i,v_i). 
    \end{align}It is easy to verify that $(U,\rho)$ is an ultrametric
   space. For each $x\in A_i$ we put
   \begin{align}\label{embmap}
    f(x):= (u_{2,0},u_{3,0}, \cdots, u_{i-1,0},f_i(x), u_{i+1,0}, u_{i+2,0},\cdots).
    \end{align}We shall prove that $f$ is a $(2+\e)$-embedding from
   the closed subset $\bigcup_{i=2}^{\infty}A_i \subseteq X$ to the ultrametric
   space $(U,\rho)$. Note that $\dim_H(\bigcup_{i=2}^{\infty}A_i)=\infty$.

   We take two arbitrary points $x\in A_i$ and $y\in A_j$
   $(i<j)$ and fix $z\in A_1$. By (\ref{ms2}) and (\ref{ms3}), we get
   \begin{align*}\dist(x,z)\leq R_i + \diam A_1 +
   \diam A_i \leq R_i+2 \leq \sqrt{1+2^{-1}\e}R_i. 
    \end{align*}Combining this inequality with (\ref{ms2}) and
   (\ref{ms3}) also implies
   \begin{align*}
    \dist(x,y)\geq \dist(y,z)-\dist(x,z)\geq R_j - \sqrt{1+2^{-1}\e}R_i \geq \frac{1}{\sqrt{1+\e}}R_j = \frac{1}{\sqrt{1+\e}}\rho(f(x),f(y)).
    \end{align*}and
   \begin{align*}
    \dist(x,y)\leq \dist(x,z)+\dist(y,z)\leq \ & \sqrt{1+2^{-1}\e}R_i
    +\sqrt{1+2^{-1}\e}R_j \\ \leq \ & \sqrt{1+\e} R_j\\
    =\ & \sqrt{1+\e} \rho(f(x),f(y)).
    \end{align*}Hence $f$ is a $(2+\e)$-embedding.

   \begin{bai}\label{case2} $\{x_i\}_{i=1}^{\infty}$ is bounded but not totally bounded.
    \end{bai}
   By taking a subsequence, we
   may assume that there exist two constants $c_1,c_2>0$ such that
   \begin{align*}
    c_1 \leq \dist (x_i,x_j)\leq c_2 \text{ for any distinct }i,j.
    \end{align*}For any $\delta>0$
  we divide $[c_1,c_2]=\bigcup_{j=1}^m I_j$ so that $\diam I_j < \delta$
  for any $j$.

  Pick $j_1\in \{1,2,\cdots, m\}$ such that $\dist(x_i, x_1)\in
  I_{j_1}$ holds for infinitely many $i$. Put
  \begin{align*}
   X_1:=\{ x_i \mid \dist(x_i,x_1) \in I_{j_1}   \}=\{x_{k_1(1)},x_{k_1(2)},\cdots \}.
   \end{align*}We then choose $j_2 \in \{1,2,\cdots,m\})$ so that
  $\dist(x_{k_1(i)},x_{k_1(1)})\in I_{j_2}$ holds for infinitely many
  $i$ and put
  \begin{align*}
   X_2:=\{ x_{k_1(i)}\in X_1 \mid \dist (x_{k_1(i)}, x_{k_1(1)})\in
   I_{j_2}\} = \{x_{k_2(1)},x_{k_2(2)}, \cdots  \}.
   \end{align*}Repeatedly we obtain a sequence $\{
  j_i\}_{i=1}^{\infty}$ whose terms are elements of the set $\{ 1,2,
   \cdots, m\}$ and $X_i=\{
  x_{k_i(1)}, x_{k_i(2)},\cdots \}$. By a pigeon hole argument we find a
  subsequence $\{ j_{h(i)}\}_{i=1}^{\infty} \subseteq \{
  j_i\}_{i=1}^{\infty}$ which is monochromatic, i.e., $j_{h(i)} \equiv
  l$ for some $l\in \{1,2,\cdots, m\}$. We then get $\dist (x_{k_{ h(i)}
  (i)}, x_{k_{h(j)}(j)})\in I_l$. Since $\diam I_l<\delta$ and
   $\lim_{i\to \infty}\diam A_{i}=0$, by choosing
  sufficiently small $\delta$ and taking a subsequence, we thereby get
   the following: There exists a number $\alpha\geq c_1$ such that 
\begin{align}\label{add1}
 \alpha \leq \dist (u,v) \leq (1+\varepsilon) \alpha \text{ for any
 }u\in A_{i} \text{ and }v\in A_{j} (i\neq j)
 \end{align}and $\diam A_{i}\leq (2+\e)^{-1}\alpha$. As in Case \ref{case1} we add an additional point $u_{i,0}$ to $f_i(A_i)$ and
   put $Y_i:=f_i(A_i)\cup \{ u_{i,0} \}$. We define the distance function
   $\tilde{\rho}_i $ on $Y_i$ as follows:
   $\tilde{\rho}_i(u,u_{i,0}):=\alpha $ for
   $u\in f_i(A_i)$ and $\tilde{\rho}_i(u,v):= \rho_i(u,v)$ for $u,v\in
   f_i(A_i)$. Since $\diam f_i(A_i)\leq (2+\e)\diam
    A_{i}\leq \alpha$, each
   $(Y_i,\tilde{\rho}_i)$ is an ultrametric space. From these $(Y_i,
   \tilde{\rho}_i)$ we construct an ultrametric space $(U,\rho)$ by
   (\ref{ultrasp}) and (\ref{ultramet}). Then a map
   $f:\bigcup_{i=2}^{\infty}A_i \to (U,\rho)$ defined by (\ref{embmap}) is a $(2+\varepsilon)$-embedding.

   \begin{bai}\label{case3} $\{x_i\}_{i=1}^{\infty}$ is totally bounded.
    \end{bai}The proof is similar to Case \ref{case1}. From the totally boundedness, by taking a subsequence, we may assume that
   $\{x_i\}_{i=1}^{\infty}$ is a Cauchy sequence. Since $\lim_{i\to
   \infty}\diam A_i =0$, the sequence $\{ A_i\}_{i=1}^{\infty}$
   Hausdorff converges to a point $x_{\infty}$. Let $\delta>0$ be
   specified later. Note that $x_{\infty}
   \notin A_i$ for any sufficiently large $i$ since $A_i$ are mutually
   disjoint closed subsets of $X$. Hence, by taking a subsequence, we
   may also assume that $\dist (A_{i},x_{\infty})/\dist (A_{i-1},x_{\infty})\leq \delta$ for each $i$. Covering $A_{i}$ by countably many
   closed subsets $\{ B_{ij}\}_{j}$ of diameter $\leq \delta \dist(A_i,x_{\infty})$ we find
   a subset $B_{ij}$ such that $\dim_H(B_{ij})\geq 2^{-1}\dim_H(A_i)$ and 
   \begin{align*}
    \frac{\diam B_{ij}}{\dist(B_{ij},x_{\infty})} \leq \frac{\diam
    B_{ij}}{\dist(A_i,x_{\infty})}\leq \delta.
    \end{align*}
   Hence by replacing $A_i$ with $B_{ij}$, we may
   assume that 
   $\diam A_{i}/\dist(A_i,x_{\infty})\leq \delta$ for every $i$.

   As in Case \ref{case1} and \ref{case2} we add an additional point $u_{i,0}$ to $f_i(A_i)$ and
   put $Y_i:=f_i(A_i)\cup \{ u_{i,0} \}$. Define the distance function
   $\tilde{\rho}_i $ on $Y_i$ by $\tilde{\rho}_i(u,u_{i,0}):=\dist(A_i, x_{\infty})$ for
   $u\in f_i(A_i)$ and $\tilde{\rho}_i(u,v):= \rho_i(u,v)$ for $u,v\in
   f_i(A_i)$. If $\delta\leq (2+\e)^{-1}$, then we have
   \begin{align*}
    \diam f_i(A_i) \leq (2+\e)\diam A_i \leq \dist(A_i,x_{\infty}),
    \end{align*}which implies that each
   $(Y_i,\tilde{\rho}_i)$ is an ultrametric space. From these
   $(Y_i,\tilde{\rho}_i)$ we define an ultrametric space $(U,\rho)$ by
   (\ref{ultrasp}) and (\ref{ultramet}). If we trace the proof of Case \ref{case1} by replacing $R_i$ with $\dist
   (A_i,x_{\infty})$, then we easily see that that a map
   $f:\bigcup_{i=2}^{\infty}A_i \to (U,\rho) $ defined by
   (\ref{embmap}) is $(2+\e)$-embedding, provided that $\delta>0$ small enough.  This completes the proof.
 \end{proof}

 \begin{proof}[Proof of Theorem \ref{finiteTHM}]Let $X:=\{x_1,x_2,\cdots
  \}$. Apply the proof of Theorem \ref{THM} by identifying each $x_i$
  with $K_i$. Note that the loss of the distortion in the proof only
  comes from (\ref{ms1}), which we can ignore in the case where
  $A_i=x_i$. Hence the space $X$ can embed into an ultrametric space
  with distortion $1+\e$. This completes the proof.
  \end{proof}

\begin{rem}\upshape After this work was completed, the author proved in \cite{Funa} that every
 proper ultrametric space isometrically embeds into $\ell_p$ for any
 $p\geq 1$. In particular a subset $S$ in Theorem
 \ref{MNDV} also embeds into $\ell_p$. Theorems \ref{finiteTHM} and \ref{THM} also hold
 in the case where the target metric space is $\ell_p$ instead of an
 ultrametric space. In fact, in the proof of Theorem \ref{THM} observe that we may assume that
 $A_i$ is compact (\cite{carleson},
  \cite[Corollary 7]{how}). Since $\bigcup_{i=2}^{\infty}A_i$ is a proper subset
 which embeds into an ultrametric space in the case of Case \ref{case1}
 and \ref{case3},
 we consider only Case \ref{case2}. Since we have
 (\ref{add1}) in Case \ref{case2} we easily see that
 $\bigcup_{i=2}^{\infty} A_i$ embeds into $\ell_p$. It was mentioned in \cite[Proposition 3.4]{Funa} that an $\ell_p$ analogue of
 Theorem \ref{hanrei} also holds.
 \end{rem}
            \begin{ack}\upshape
        The author would like to express his thanks to Mr. Takumi Yokota for
             his suggestion regarding the two theorems in this paper and
             Mr. Ryokichi Tanaka for discussion. The author also thanks
             Professor Manor Mendel and Professor Assaf Naor for their
             useful comments. The author is also indebted to an anonymous referee for carefully reading
             this paper and helpful comments. 
 \end{ack}

	\end{document}